\documentclass[10pt]{article}
\usepackage[francais,english]{babel}
\usepackage{amsmath}
\usepackage{amsfonts}
\usepackage{amssymb}
\usepackage{amsthm}
\usepackage{graphics}
\usepackage{amscd}
\usepackage{epsfig}
\usepackage{color}
\newcommand{\Z}{\mathbb{Z}}
\newcommand{\R}{\mathbb{R}}

\newcommand{\C}{\mathbb{C}}

\newcommand{\U}{\mathbb{U}}
\newcommand{\E}{\mathbb{E}}

\newcommand{\ind}{{\bf 1}}
\renewcommand{\P}{\mathbb{P}}

\DeclareMathOperator{\Id}{Id}

\DeclareMathOperator{\SLE}{SLE}

\title{Critical percolation in annuli and $\SLE_6$}
\author{Julien Dub\'edat\footnote {Universit\'e Paris-Sud}}

\newtheorem{Conj}{Conjecture}
\begin{document}
\maketitle

\begin{abstract}
Building on the identification of the scaling limit of the critical percolation
exploration process as a Schramm-Loewner Evolution,
we derive a PDE characterization for the crossing probability of an annulus.
\end{abstract}

\section {Introduction}

Percolation is arguably the simplest example of a planar ``critical'' model,
i.e. a random planar graph (generally speaking a subgraph of a regular
lattice, eventually with spins on sites) with a small mesh such that the
probabilities of geometrically meaningful macroscopic events have non
degenerate limits when the mesh goes to zero. Recall that percolation consists in
removing each edge (or each vertex) in a lattice with a given
probability $p$; bond percolation on the square lattice and site
percolation on the triangular lattice are critical for $p=1/2$.

The probabilities of macroscopic events are generally believed to be conformally invariant: 
the limiting probability when the mesh of the lattice goes to zero should
not change if one applies a conformal equivalence to the corresponding
geometric configuration. An important related example is the 
conformal invariance of hitting probabilities (harmonic measure) for
planar Brownian motion, which is the scaling limit of simple random walks.
Using techniques from Conformal Field Theory (CFT) and Coulomb gas models,
physicists have proposed several intriguing formulas for the limiting probabilities of
critical percolation. Unfortunately, it does not seem easy to make their arguments rigorous. 
We now review some of these predictions.

\begin{itemize}
\item Simply connected domains\\
Mark four points on the boundary of a simply connected domain to get a
conformal rectangle; up to conformal equivalence, there is only one
degree of freedom, namely the aspect ratio of the rectangle. The
probability that two opposite sides of the rectangle are connected by
a percolation cluster is given by Cardy's formula (\cite{Ca}). Watts' formula
(\cite{Wa}) describes the probability of a double crossing, i.e. the two pairs of opposite sides
are all connected by one cluster. Cardy has also studied the occurrence
of $n$ disjoint clusters connecting two opposite sides (\cite{Ca1}).
\item Multiply connected domains\\
A $n$-connected domain is a plane connected domain the complement in
$\C$ of which has $n$ connected components (its fundamental
group is the free group on $(n-1)$ generators). A natural question
concerns the probability that two given connected components of the
boundary are connected by a cluster inside the domain. The simplest
case is the existence of a crossing from the outer to the inner
boundary of an annulus. More generally, a law for the number of
crossings of alternate colors has been proposed by Cardy (\cite{Ca2}). 
\item Compact Riemann surfaces\\
Here there is no boundary; tori are the simplest example. Since the
conformal structure is of crucial importance, it may be better to
think of (complex) elliptic curves. Pinson (\cite{Pi}) has studied the image of
the (first) homology group of clusters of a given color in the homology
group of the torus (i.e. a random subgroup of $\Z\times\Z)$.
\end{itemize}

Since a continuous random graph does not really make sense, one has first to 
clarify what these limiting probabilities do correspond to. 
A natural way to understand the ``continuous scaling limit of critical percolation''
is to focus on the interfaces between clusters, which are random curves.
The limiting probabilistic objects (random curves) are the Schramm-Loewner Evolution (SLE),
introduced by Schramm in \cite{S0} (see \cite{RS01,W1} for some background
on SLE). It describes the only possible conformally invariant scaling limits
of these interfaces. 

For critical site percolation on the triangular lattice, Smirnov
(\cite{Sm1}) has proved that the cluster interface indeed converges
to the so-called chordal $\SLE_6$ process. His proof in fact relies on
establishing directy Cardy's formula for single crossings in conformal
rectangles. 

Hence it may be interesting to derive some percolation probabilities
in the $\SLE_6$ framework. Cardy's formula itself is easily derived
for the $\SLE_6$, as pointed out by Schramm (see for instance \cite {W1}).
The general idea is that if one wants to compute probabilities of 
macroscopic events in terms of the $\SLE_6$ process, one uses the fact that 
the conditional probabilities (using the filtration associated to the $\SLE_6$)
of this event are a martingale. This leads (usually) to a partial differential equation
(in terms of the ``conformal coordinates'' of the problem), that characterize
fully these probabilities. This method has been used in 
\cite {S1} to derive a new formula for critical percolation, and also
in \cite {LSW1,LSW4,LSW6,LSW7} to derive the value of the corresponding 
critical exponents (that describe the asymptotic decay of the probabilities of 
some events).

In the present paper we derive analytic characterizations of some percolation
probabilities in annuli using Smirnov's result on convergence of the
percolation exploration process to $\SLE_6$. The computations
are somewhat tedious, but the underlying probabilistic ideas are extremely
simple, and may be of use for more general problems (in particular
for $n$-connected domains, $n\geq 3$). 

Let $U$ be a $n$-connected domains and $\partial$ be one connected
component of its boundary (a Jordan closed curve, say). Pick two
points $x$ and $y$ on $\partial$ and set the following boundary
conditions on $\partial$: the arc $(x,y)$ is colored in blue and the
arc $(y,x)$ in yellow (in clockwise order). A {\it chordal event} is a
percolation event depending only on $(U,x,y)$. For instance, consider
the event that there exists a blue crossing between $(x,y)$ and
$\partial_1$ and a yellow crossing between $(y,x)$ and $\partial_2$,
where $\partial_i$ are connected components of $\partial U$. Then
``grow'' a small percolation (resp. $\SLE_6$) hull at $x$. Let $K_t$ be
this hull ($t$ is a measure of its size) and $x_t$ be its ``tip'' 
(see e.g. \cite{Sm1}); then
$(U\backslash K_t,x_t,y)$ is the perturbed domain. For a well chosen
chordal event, the event holds in $(U,x,y)$ if and only if it holds in
$(U\backslash K_t,x_t,y)$.

Now consider the set of $(g+1)$-connected domain with two points marked on
a connected component of the boundary modulo conformal
equivalence (a Teichm\"uller space); it is a manifold of dimension
$(3g-1)$ if $g\geq 2$. Then $t\mapsto(U\backslash K_t,x_t,y)$ (modulo conformal
equivalence) should be, up to time change, a diffusion in the
Teichm\"uller space; as percolation is local, it does not ``feel'' the
boundary before actually touching it. The probability of the chordal event defines a
function on this space and this function is harmonic for the diffusion. If one is
able to compute a SDE for the diffusion and, crucially, to work out
boundary conditions for the harmonic function, this yields an analytic
characterization of the chordal event probability (as a function on
the Teichm\"uller space).

In the paper, we essentially carry out this program for annuli. In the
first section, we recall some facts of complex analysis on the annulus
(mainly the solution to the Dirichlet problem). Next, we derive the
SDE for the diffusion in the Teichm\"uller space (which is isomorphic
to $\R_-\times (0,2\pi)$). In the third section, we characterize the
law of the number of nested clusters of alternate colors wrapped
around the inner disk of an annulus, and make the connection with
Cardy's results (see \cite{Ca2}).

Defining SLE on Riemann surfaces has also been recently (and independently) 
the subject of \cite {DZ,FK}, in different settings.

{\bf Acknowledgments.} I wish to thank Wendelin Werner for his help
and advice, and Vincent Beffara for very enlightening discussions.

\section {Annuli}
Let $q<1$; define the annulus
$$A_q=\{z: q<|z|<1\}.$$
It is classical that for $q\neq q'$, the annuli $A_q$ and $A_{q'}$ are
not conformally equivalent, and that the conformal automorphisms of
$A_q$ are the maps $z\mapsto uz$ and $z\mapsto quz^{-1}$, $u\in
\U$. Moreover, any doubly connected domain (i.e. a connected domain the
complement in $\C$ of which has two connected components) is
conformally equivalent to an annulus. Thus, one may identify the
Teichm\"uller space of doubly connected domains with the set
$\{A_q\}_{0<q<1}$.

Let us briefly recall Villat's solution of the Dirichlet problem in an annulus
\cite{Vil12}. Let $\phi$, $\psi$ be two continuous, real-valued
$2\pi$-periodic functions. The Dirichlet problem consists in finding a
real harmonic function $f$ defined in $A_q$ with boundary values given
by $\phi$ and $\psi$. The classical Dirichlet problem in the disk may
be solved using the Poisson kernel. In the case of annuli, one may also exhibit
a kernel, which involves
elliptic functions (see e.g. \cite{Cha}). Let $\omega_1$, $\omega_2$
be two numbers such that $\omega_1$ is real positive, $\omega_2$ is
imaginary, and 
$$q=\exp(-\frac{\pi\omega_2}{i\omega_1})$$
We shall consider elliptic functions with basic periods
$(2\omega_1,2\omega_2)$. Recall that the Weierstrass zeta-function is
given by:
$$\zeta(z)=\zeta(z;2\omega_1,2\omega_2)=\frac 1z+\sum_{\omega\neq
  0}\left(\frac{1}{z-\omega}+\frac 1\omega+\frac z{\omega^2}\right)$$
where the sum is on the vertices of the lattice
$2\omega_1\Z+2\omega_2\Z$. Then, for any $z$,
$\zeta(z+2\omega_1)=\zeta(z)+2\eta_1$,
$\zeta(z+2\omega_2)=\zeta(z)+2\omega_2$, with
$\eta_1=\zeta(\omega_1)$, $\eta_2=\zeta(\omega_2)$. Let
$\omega_3=\omega_1+\omega_2$ and $\eta_3=\eta_1+\eta_2$. Define also:
$\zeta_3(z)=\eta_3-\zeta(\omega_3-u)$.

Note that $z\mapsto \log|z|$ defines a real harmonic function on an
annulus which is not the real part of an holomorphic function on this
annulus. So assume that:
$$\int_0^{2\pi}\phi(\theta)d\theta=\int_0^{2\pi}\psi(\theta)d\theta$$
Then one may define
\begin {eqnarray*}
\Omega(z)&=&\frac{i\omega_1}{\pi^2}\int_0^{2\pi}\phi(\theta)\zeta\left(\frac{\omega_1}{i\pi}\log(z)-\frac{\omega_1}{\pi^2}\theta\right)d\theta\\
&&
-\frac{i\omega_1}{\pi^2}\int_0^{2\pi}\psi(\theta)\zeta_3\left(\frac{\omega_1}{i\pi}\log(z)-\frac{\omega_1}\pi\theta\right)d\theta.
\end {eqnarray*}
Along a clockwise loop in $A_q$, the first integral is increased by
$i\omega_1/{\pi^2}\int_0^{2\pi}\phi(\theta)\eta_1d\theta$ and the second
integral is increased by the same amount; so there is no problem with
the logarithms determination.

This function $\Omega$ is holomorphic on $A_q$; moreover, for all $\theta$, 
$$\Re\Omega(\exp(i\theta))=\phi(\theta)\text{\ \ and\ \ } \Re
\Omega(q\exp(i\theta))=\psi(\theta)$$

We apply this result with the conditions
$\phi(\theta)d\theta=2\pi\delta_{\theta_0}$, and $\psi(\theta)=1$ for
all $\theta$, thus getting a holomorphic function $\Omega$ which is
well defined up to the addition of an imaginary constant. Let $x=\exp(i\theta_0)$ and
$y\in\U$, $y\neq x$. The holomorphic vector field: 
\begin{align*}
V_{x,y}(z)=&\frac{2i\omega_1}\pi
z\left(\zeta\left(\frac{\omega_1}{i\pi}\log(z/x)\right)-\zeta\left(\frac{\omega_1}{i\pi}\log(y/x)\right)\right)\\
&-\frac{2i\omega_1}\pi z\left(\int_0^{2\pi}\zeta_3\left(\frac{\omega_1}{i\pi}\log(ze^{-i\theta})\right)-\zeta_3\left(\frac{\omega_1}{i\pi}\log(ye^{-i\theta})\right)\frac{d\theta}{2\pi}\right)
\end{align*}
is equal to $z\Omega(z)$, so it is well defined by the following properties:
\begin{itemize}
\item $V_{x,y}$ is holomorphic on $A_q$,
\item $V_{x,y}$ may be extended continuously to the boundary, except
  at $x$,
\item $\Re(V_{x,y}(z)/z)$ is constant on $q\U$,
\item $\Re(V_{x,y}(z)/z)=0$ on $\U\backslash\{x\}$, 
\item $V_{x,y}(y)=0$, and
\item $V_{x,y}$ has a simple pole at $x$ with residue $-2x^2$.
\end{itemize}

\section{$\SLE_6$ in an annulus}

A chordal $\SLE_\kappa$ going from $x\in\U$ to $1$ in the unit disk $D$
may be defined (up to linear time change) by the equations:
$$dg_t(w)=L_{W_t}(g_t(w))dt$$
where $W$ is $\sqrt\kappa$ times a standard real Brownian motion
starting from $i(x+1)/(1-x)$ and the holomorphic vector field $L_y$ is
defined on $D$ by:
$$L_y(w)=-\frac{i(1-w)^2}{i(1+w)/(1-w)-y}$$
The map $g_t$ is a conformal equivalence between $D\backslash K_t$ and
$D$ that fixes $1$, where $K_t$ is the $\SLE$ hull at time $t$. The
time parameter corresponds to half of the capacity of the image of
$K_t$ under the homography $w\mapsto i(w+1)/(1-w)$ seen from infinity
in the half-plane.

Let $q<1$. For small enough $t$, $D\backslash g_t(qD)$ is a
doubly connected set; hence there exists a unique $q_t\ge q$ and a unique conformal
equivalence $h_t$ between $D\backslash g_t(qD)$ and $A_{q_t}$ such that
$h_t(1)=1$. Let $\xi_t=(W_t-i)/(W_t+i)$, $\lambda_t=h_t(\xi_t)$,
and $f_t=h_t\circ g_t$.

Now $w\mapsto (df_t)(f_t^{-1}(w))$ defines a
holomorphic vector field on $A_{q_t}$. It is easily seen that this
vector field is proportional to $V_{\lambda_t,1}$ (we omit the $q_t$ parameter). Thus, there exists
a function $a_t$ such that:
$$d(f_t(w))=V_{\lambda_t,1}(f_t(w))da_t$$
Moreover, $d\log(q_t)=da_t$, so one may pick $a(t)=\log(q(t))$. Using
the chain rule, one gets:
$$d(f_t(w))=(dh_t)(g_t(w))+h'_t(g_t(w))dg_t(w)$$
It follows that:
$$dh_t(w)=V_{\lambda_t,1}(h_t(w))da_t-L_{W_t}(w)h'_t(w)$$
As the poles in this expression should cancel out, necessarily:
$$da_t=h'_t(\xi_t)^2\frac{(1-\xi_t)^4}{4\lambda_t^2}$$
and this is indeed real. Now let $g_t(w)=\xi_t+u$, with small $u$. Then,
$$
L_{W_t}(\xi_t+u)=-\frac{(1-\xi_t)^4}{2u}+\frac 32(1-\xi_t)^3+O(u),$$
\begin {eqnarray*}
\lefteqn {h'_t(\xi_t+u)L_{W_t}(\xi_t+u)} \\
&=&
-h'_t(\xi_t)\frac{(1-\xi_t)^4}{2u}-h''_t(\xi_t)\frac{(1-\xi_t)^4}{2}+\frac
32 h'_t(\xi_t)(1-\xi_t)^3+O(u)
\end {eqnarray*}
and 
\begin {eqnarray*}
\lefteqn {\zeta\left(\frac{\omega_1}{i\pi}\log\left(\frac {h_t(\xi_t+u)}{\lambda_t} \right)\right)}\\
&=&
\frac{i\pi\lambda_t}{\omega_1h'_t(\xi_t)}\left(\frac
  1u-\frac{\lambda_t}{2h'_t(\xi_t)}\left(\frac {h''_t(\xi_t)}{\lambda_t}-\frac{h'_t(\xi_t)^2}{\lambda_t^2}\right)\right)+O(u) 
.\end {eqnarray*}
Finally
\begin{eqnarray*}
\lefteqn{V_{\lambda_t,1}(h_t(\xi_t+u))} \\
&=&\lambda_t\left[-\frac{2\lambda_t}{h'_t(\xi_t)}\left(\frac
  1u-\frac{\lambda_t}{2h'_t(\xi_t)}\left(\frac    {h''_t(\xi_t)}{\lambda_t}-\frac{h'_t(\xi_t)^2}{\lambda_t^2}\right)\right)+\frac{2i\omega_1}\pi\zeta\left(\frac{\omega_1}{i\pi}\log(\lambda_t)\right)\right.\\
&&+\left.\frac{2i\omega_1}\pi\int_0^{2\pi}\left(\zeta_3\left(\frac{\omega_1}{i\pi}\log
     ( e^{-i\theta})\right)-\zeta_3\left(\frac{\omega_1}{i\pi}\log
     (\lambda_t  e^{-i\theta})\right)\right)\frac{d\theta}{2\pi}\right]\\
&&-2\lambda_t+O(u)
\end{eqnarray*}
Notice that:
$$\int_0^{2\pi}\left(\zeta_3\left(\frac{\omega_1}{i\pi}\log
     ( e^{-i\theta})\right)-\zeta_3\left(\frac{\omega_1}{i\pi}\log
     (\lambda_t  e^{-i\theta})\right)\right)\frac{d\theta}{2\pi}=-\frac{\log(\lambda_t)}{2i\pi}2\eta_1$$
Hence $w\mapsto dh_t(w)$ is smooth at $\xi_t$ and:
\begin{align*}
(dh_t)(\xi_t)=&\frac{h'_t(\xi_t)^2(1-\xi_t)^4}{4\lambda_t}\left[
\frac{2i\omega_1}{\pi}\zeta\left(\frac{\omega_1}{i\pi}\log(\lambda_t)\right)-\frac{2\omega_1\log(\lambda_t)\eta_1}{\pi^2}\right]\\
&+\frac 34 h''_t(\xi_t)(1-\xi_t)^4-\frac
{3h'_t(\xi_t)^2(1-\xi_t)^4}{4\lambda_t}-\frac 32 h'_t(\xi_t)(1-\xi_t)^3
\end{align*}
Recall that $\xi_t=(W_t-i)/(W_t+i)$, $dW_t=\sqrt\kappa dB_t$, where $B$
is a standard real Brownian motion. Hence:
$$d\xi_t=\frac{(1-\xi_t)^2}{2i}\sqrt\kappa
dB_t+\frac{(1-\xi_t)^3}4\kappa dt$$
Since $\lambda_t=h_t(\xi_t)$, an appropriate version of It\^o's
formula yields:
\begin{eqnarray*}
d\lambda_t &=&(dh_t)(\xi_t)+h'_t(\xi_t)d\xi_t+\frac
12 h''_t(\xi_t)d\langle\xi_t\rangle\\
&=&\frac{h'_t(\xi_t)(1-\xi_t)^2}{2i}\sqrt\kappa dB_t
\\
&&+\frac{h'_t(\xi_t)^2(1-\xi_t)^4}{4\lambda_t}\left[
\frac{2i\omega_1}{\pi}\zeta\left(\frac{\omega_1}{i\pi}\log(\lambda_t)\right)-\frac{2\omega_1\log(\lambda_t)\eta_1}{\pi^2}\right]dt\\
&&
+\left(h''_t(\xi_t)\frac{(1-\xi_t)^4}2-h'_t(\xi_t)(1-\xi_t)^3\right)\left(\frac
32-\frac\kappa 4\right)dt
\\
&&
-\frac{3h'_t(\xi_t)^2(1-\xi_t)^4}{4\lambda_t}dt
\end{eqnarray*}
Let $\kappa=6$ and $\exp(i\nu_.)=\lambda_.$, $\nu_.\in [0,2\pi]$. We perform a time change
using the increasing function $t\mapsto a(t)$:
$$d\nu_a=-\sqrt\kappa d\tilde B_a+\frac
{2\omega_1}{\pi}\left(\zeta(\omega_1\frac{\nu_a}{\pi})-\frac{\nu_a}{\pi}\zeta(\omega_1)\right)da$$
where $(\tilde B_a)_{a\geq 0}$ is a standard real Brownian
motion. Note that the zeta-functions, as well as the half-period $\omega_1$
depend implicitly on $a$. Using Jacobi's theta-function, one may
rewrite this SDE as (see e.g. \cite{Cha}, Theorem V.2):
$$d\nu_a=-\sqrt\kappa d\tilde B_a+\frac 1\pi\frac{\theta'}\theta\left(\frac{\nu_a}{2\pi},-\frac{ia}{\pi}\right)da$$
where $\theta'$ designates the derivative of the bivariate function
$\theta(v,z)$ with respect to the first variable. Note that this SDE
is invariant under $\nu\leftrightarrow 2\pi-\nu$, which is
obvious from the definition of $\nu$.

The principle of this computation is closely related to the approach
of locality/restriction in \cite{LSW3}. The fact that one gets a
(time-inhomogeneous) diffusion for $\nu$ is a feature of locality for
$\kappa=6$. Note also that we could have begun with any configuration
conformally equivalent to $(A_q,\lambda,1)$, getting the same dynamics
for $\nu$, hence the same percolation probabilities. 

\section{Crossing of an annulus}

The previous computations may be used to study various critical percolation
probabilities in the annulus. For instance, consider the following
crossing probability: $F(\nu,a)$ is the probability that there exists a blue
crossing between the arc $(1,\exp(i\nu))\subset\U$ and the inner circle
$q\U$ in the annulus $A_q$ and a yellow crossing between
$(\exp(i\nu),1)$ and $q\U$, with $a=\log(q)$.  It
follows from the previous section that:
\begin{equation}\label{PDE}
3F_{\nu,\nu}+\frac
1\pi\frac{\theta'}\theta\left(\frac{\nu_a}{2\pi},-\frac{ia}{\pi}\right)F_\nu+F_a=0
\end{equation}
and $F$ satisfies the boundary conditions: $F(0,a)=F(2\pi,a)=0$, $F(\nu,0)=1$
for all $a<0$, $\nu\in (0,2\pi)$. This last condition is a consequence
of the Russo-Seymour-Welsh theory.

In the rest of this section, we discuss events related to nested
circuits of alternate colors around the inner disk of an annulus for
critical percolation. 

Consider $(\nu_a,a)_a$ as a two dimensional diffusion in the
half-strip $S=\{z:\Re z<0, 0<\Im z<2\pi\}$. This half-strip may be
identified with the Teichm\"uller space of doubly connected plane
domains with two marked boundary points (on the same connected
component of the boundary) modulo conformal equivalence. The diffusion
is stopped when it hits $\R_-\cup [0,2i\pi]$ and is instantaneously
reflected on $2i\pi+\R_-$. For any $z\in S$, this defines a harmonic
measure seen from $z$ on $\R_-\cup [0,2i\pi]$. The restriction of this
probability measure to
$\R_-$ is a finite measure (not a probability measure), which we shall
note $\tilde K(z,.)$. For $y\in\R_-$, note $K(y,.)=\tilde
K(y+2i\pi,.)$. Thus we have defined a defective Markov kernel on $\R_-$
(identified with doubly connected domains modulo conformal
equivalence). From a probabilistic point of view, this is equivalent
to a (discrete time) Markov chain on $\R_-\cup\{\partial\}$, where
$\partial$ is a cemetery state. The number of steps taken by this
Markov chain (starting from $y<0$) before reaching $\partial$ corresponds to the number of
downcrossings and upcrossings of $\nu$ (starting from $2\pi$ at time $y$). For
a given $y$, $(\nu,a)\mapsto \tilde K(a+i\nu,y)$ is a solution of the PDE (\ref{PDE})
with boundary conditions: $\tilde K(i\nu,y)=0$ for all $\nu\in (0,2\pi)$,
and, on $\R_-$, $\tilde K(.,y)=\delta_y$, the Dirac mass at $y$. On
$2i\pi+\R_-$, there is a Neumann boundary condition:
$\tilde K_\nu(2i\pi+x,y)$ for all $x<0$, just as in \cite{LSW7}, Lemma 2.3.

We now interpret downcrossings and upcrossings of $\nu$ in terms of
the exploration process. See also \cite{LSW7} for a discussion of
percolation events in annuli in relation with the exploration process (and nice
figures !). For the sake of simplicity, consider critical site percolation on the
triangular lattice: each vertex of the triangular lattice (or each
hexagon of a honeycomb lattice) is colored in blue or yellow with
probability $1/2$. Consider a portion of this lattice (with small
mesh) that approximates the annulus $A_q$. Two points are marked on
the outer boundary, say $x$ and $1$; the arc $(1,x)$ is colored in
blue and the arc $(x,1)$ is colored in yellow. The exploration process
from $x$ to $1$ (which is well defined as long as it does not hit
$q\U$) is the path with only blue hexagons on its left-hand and yellow ones on
its right-hand. The exploration process completes a clockwise loop if
there is a blue circuit of hexagons around the inner disk $qD$ and a
counterclockwise loop if there is a yellow circuit. It hits the inner
disk before completing a circuit if there is a blue crossing from $(1,x)$
to $q\U$ and a yellow one from $(x,1)$ to $q\U$. In the following,
a circuit will always be a circuit around $qD$, and a crossing will
always be a crossing between $\U$ and $q\U$ in $A_q$.

In the continuous setting, starting from $\exp(i\nu_a)=x$, if $\nu_.$
reaches $0$, then the boundary ``seen from the inner disk'' is
completely yellow; this means that the exploration process has
completed a counterclockwise loop (including the yellow arc
$(x,1)$). Hence, starting from $\nu_a=2i\pi$, the process $\nu$ reaches
$0$ before time $0$ if and only if there is a yellow circuit inside the
annulus (the outer circle is blue). Then $\nu$ reflects
instantaneously on $\R_-$, which means that the exploration process
proceeds towards $qD$ in a conformal annulus the outer boundary of
which (i.e. the
circuit) is yellow. Thus a downcrossing for $nu$ means that there exists a
yellow circuit (here $\U$ is blue); a downcrossing followed by an
upcrossing means that there exists a yellow circuit and a blue circuit
nested in the yellow one; and so on.

Now set free boundary conditions. Let $N$ denote the event that there is
no circuit (blue or yellow); equivalently, there is a blue crossing
and a yellow one. For $n\geq 1$, let $B_n$ be the event that there is exactly $n$
disjoint clusters of alternate colors wrapped around $qD$ and the
outermost is blue; $Y_n$ is the corresponding events with colors
changed (see figure). 

\begin{figure}[htbp]
\begin{center}
\input{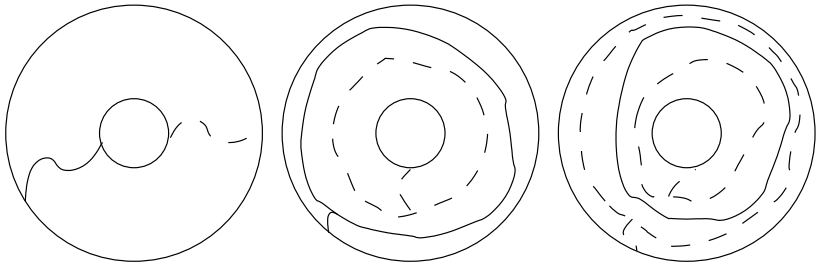}
\end{center}
\caption{The events $N$, $B_2$, $Y_3$} 
\end{figure}

Obviously $\P(B_i)=\P(Y_i)$ and:
$$\P(N)+2\sum_{n=1}^\infty\P(B_i)=1$$
It follows from the previous discussion that the probability $c_n$ that
$\nu$ starting from $y+2i\pi$ makes at least $n$ alternate downcrossings and
upcrossings is:
$$c_n=(K^{\star n})(y,\ind)$$
where $\star$ designates the convolution of Markov kernels and $\ind$
is the constant function. Note that we have analytically characterized
$K$. Here one should be cautious because of the boundary
conditions. The number $c_n$ is the probability to see at least $n$
circuits of alternate colors, not counting the outermost one if it is
blue. This outermost circuit is either blue or yellow, so:
\begin{align*}
c_n&=\sum_{k=n}^\infty\P(Y_k)+\sum_{k=n+1}^\infty\P(B_k)\\
&=\P(B_n)+2\sum_{k=n+1}^\infty\P(B_k)
\end{align*}
Determining $\P(N),\P(B_i)$ from the $c_i$ is then a trivial problem
of (infinite dimensional) linear algebra. In the Fr\'echet space
of rapidly decaying series, consider the vectors $v=(\P(N),\P(B_1),\dots,\P(B_n)\dots)$
and $w=(1,c_1,\dots,c_n,\dots)$, and the bounded operator:
$$M=\Id+2J+2J^2+\dots+2J^n+\dots$$
where $J$ is the left shift: $J(u_0\dots u_n\dots)=(u_1\dots
u_{n+1}\dots)$. Then $Mv=w$. As (a little formally)
$M=2(\Id-J)^{-1}-\Id$, one gets $M^{-1}=\Id-2(\Id+J)^{-1}$. Hence:
$$\P(B_i)=c_i+2\sum_{n=1}^\infty(-1)^nc_{n+i}$$ 
with the conventions $B_0=N$, $c_0=1$. Note that $(c_n)$ decays
exponentially: once a circuit has been completed, the ``new annulus'' has
a greater modulus that the initial one, and the greater the modulus,
the lesser the probability that there exists circuits. Hence $c_n(q)<c_1(q)^n$, and
$c_1(q)<1$ from the RSW theory (one may also argue using the BK
inequality rather than the previous ``renewal'' argument). Similarly,
$(\P(B_n))$ decays exponentially, which justifies the previous formal argument.

Let $M$ be the number of downcrossings and upcrossings completed by
the diffusion $\nu$ starting from $2i\pi+\log q$, and $\epsilon$ be
the random sign:
$$\epsilon=1+2\sum_{i=1}^\infty \ind_{\{M\geq i\}}$$ 
The sum is finite a.s. From the preceding discussion,
$\P(N)=\E(\epsilon)$ as functions of $a=\log q$. Note that
$\epsilon=1$ if the last visited boundary half-line is the top one, and
$\epsilon=-1$ in the other case. Since $\nu$ and $2i\pi-\nu$ have the
same law, if follows that:
$$\E_{i\pi+\log q'}(\epsilon)=0$$  
for all $q'\in (0,1)$. Moreover, if $\nu$ does not reach the half-line
$i\pi+\R_-$, then $\epsilon=1$ (since paths are continuous). Hence, by the Markov property:
$$\P(N)(q)=\E_{2i\pi+\log q}(\epsilon)=\P_{2i\pi+\log q}(\nu\text{\ does not reach }i\pi+\R_-)$$
So we can characterize $\P(N)$ as the solution of a first exit
problem: consider the diffusion $\nu$ starting from $z$, $\Re z<0$,
$0<\Im z<\pi$ (using
the symmetry $\nu\leftrightarrow 2i\pi-\nu$), stopped at time $\tau$
when it exits the half-strip $\{z:\Re z<0, 0\leq\Im z<\pi\}$, the bottom part
of the boundary being instantaneously reflecting; either $\Im
\nu_\tau=\pi$ or $\Re \nu_\tau=0$ a.s. Then:
$$\P(N)(q)=\P_{\log q}(\Re\nu_\tau=0)$$ 
Let $H(z)=\P_z(\Re\nu_\tau=0)$, $z=a+i\nu$, defined on the half-strip $\{z:\Re
z<0, 0<\Im z<\pi\}$. It is easily seen that in the interior of the
domain,
$$3H_{\nu,\nu}+\frac
1\pi\frac{\theta'}\theta\left(\frac{\nu_a}{2\pi},-\frac{ia}{\pi}\right)H_\nu+H_a=0 
$$ 
and $H$ satisfies the boundary conditions: $H$ equals 1 on
$i(0,\pi)$, 0 on $i\pi+\R_-$, and the normal derivative of $H$
vanishes on $\R_-$.This mixed
Dirichlet-Neumann problem characterizes completely the function
$q\mapsto \P(N)(q)$. This should be compared with the
following formula, derived by Cardy (\cite{Ca2}) using Coulomb gas
techniques, in link with Conformal Field Theory:

\begin{Conj}[Cardy] 
Let $\tau=-ia/2\pi$. With the previous notations:
$$\P(N)=\sqrt 3\frac{\eta(\tau)\eta(6\tau)^2}{\eta(3\tau)\eta(2\tau)^2}$$
where $\eta$ designates Dedekind's eta function.
\end{Conj}

So far, we have not (yet?) been able to derive this from our
characterization of $\P(N)$. In a subsequent paper (\cite{Dub}), we intend to study
percolation problems in multiply-connected domains, as well as
$\SLE_{8/3}$ in such domains.

-----------------------

Laboratoire de Math\'ematiques, B\^at. 425

Universit\'e Paris-Sud, F-91405 Orsay cedex, France

julien.dubedat@math.u-psud.fr


\begin{thebibliography}{99}

\bibitem[Ahl79]{Ahl} L. Ahlfors, {\it Complex Analysis}, 3rd edition,
  McGraw-Hill, 1979
\bibitem[Ca92]{Ca}
{J.L. Cardy, {\it Critical percolation in finite geometries},
J. Phys. A, 25, L201--L206, 1992}
\bibitem[Ca01]{Ca1} J.L. Cardy, {\it Conformal invariance and
    percolation}, preprint, arXiv:math-ph/0103018, 2001
\bibitem[Ca02]{Ca2} J.L. Cardy, {\it Crossing Formulae for Critical
    Percolation in an Annulus}, preprint, arXiv:math-ph/0208019v4, 2002
\bibitem[Cha84]{Cha} K. Chandrasekharan, {\it Elliptic functions},
  Grundlehren der mathematischen Wissenschaften 281, Springer-Verlag,
  1984
\bibitem[Dub03]{Dub} J. Dub\'edat, in preparation
\bibitem[FK03]{FK}{R. Friedrich, J. Kalkkinen, in preparation}
\bibitem[Gri97]{G1} G.R. Grimmett, {\it Percolation and disordered
    systems}, in {\it Lectures on Probability Theory and Statistics,
    Ecole d'\'et\'e de probabilit\'es de Saint-Flour XXVI}, Lecture
  Notes in Mathematics 1665, Springer-Verlag, 1997
\bibitem[Law01]{L1} G. Lawler, {\it An Introduction to the Stochastic
    Loewner Evolution}, preprint, 2001
\bibitem[LSW01a]{LSW1} G. Lawler, O. Schramm, W. Werner, {\it
    Values of Brownian intersection exponents I: Half-plane
    exponents}, Acta Math. 187, 237--273, 2001
\bibitem[LSW01b]{LSW4} G. Lawler, O. Schramm, W. Werner, {\it
    Values of Brownian intersection exponents II: Plane
    exponents}, Acta Math. 187, 275--308, 2001
\bibitem[LSW02a]{LSW6} G. Lawler, O. Schramm, W. Werner, {\it
    Values of Brownian intersection exponents III: Two-sided
    exponents}, Ann. Inst. H. Poincar\'e Probab. Statist., 38, no 1,
  109-123, 2002
\bibitem[LSW02b]{LSW2} G. Lawler, O. Schramm, W. Werner, {\it
    Conformal Invariance of planar loop-erased random walks and
    uniform spanning trees}, arXiv:math.PR/0112234, Ann. Probab., to appear
\bibitem[LSW02c]{LSW5} G. Lawler, O. Schramm, W. Werner,
{\it On the scaling limit of planar self-avoiding walk,}
 math.PR/0204277, 2002 
\bibitem[LSW02d]{LSW3} G. Lawler, O. Schramm, W. Werner, {\it
    Conformal restriction. The chordal case}, arXiv:math.PR/0209343,
    J. Amer. Math. Soc., to appear
\bibitem[LSW02e]{LSW7}  G. Lawler, O. Schramm, W. Werner, {\it One-arm
    exponent for critical 2D percolation}, Electr. J. Probab. 7, no 2,
  2002
\bibitem[Pin94]{Pi} H. Pinson, {\it Critical percolation on the
    torus}, J. Statist. Phys., 75, no 5-6, pp 1167--1177, 1994
\bibitem[RevYor94]{RY} D. Revuz, M. Yor, {\it Continuous martingales and
    Brownian motion}, 2nd edition, Grundlehren der mathematischen wissenschaften
  293, Springer Verlag, 1994
\bibitem[RohSch01]{RS01} S. Rohde, O. Schramm, {\it Basic Properties
    of SLE}, arXiv:math.PR/0106036, 2001
\bibitem[Sch00]{S0} O. Schramm, {\it Scaling limits of loop-erased random walks and
    uniform spanning trees}, Israel J. Math., 118, 221--288, 2000 
\bibitem[Sch01]{S1} O. Schramm, {\it A percolation formula},
  Electr. Comm. Probab. 6, 115--120.  
\bibitem[Smi01]{Sm1} S. Smirnov, {\it Critical percolation in the
    plane. I. Conformal Invariance and Cardy's formula II. Continuum
    scaling limit}, in preparation, 2001
\bibitem[SW01]{SW} S. Smirnov, W. Werner, {\it Critical exponents for
    two-dimensional percolation}, Math. Res. Lett., 8, pp 729--744, 2001
\bibitem[Vil12]{Vil12} H. Villat, {\it Le probl\`eme de Dirichlet dans
    une aire annulaire}, Rend. circ. mat. Palermo, pp 134--175, 1912
\bibitem [Wa96]{Wa} G.M.T. Watts, {\it A crossing probability for
    critical percolation in two dimensions}, J. Phys. A:
  Math. Gen. 29, pp 363--368, 1996
\bibitem[Wer02]{W1} W. Werner, {\it Random planar curves
    and Schramm-Loewner evolution}, Lecture Notes of the 2002 St-Floor
  summer school, Springer, to appear
\bibitem[Z03]{DZ} {D. Zhan, in preparation.}
\end{thebibliography}
\end{document}